\documentclass[final,leqno,onefignum,onetabnum]{amsart}
\usepackage{amssymb,amsmath,amstext}
\usepackage{latexsym}
\title{On efficiently computing the eigenvalues of limited-memory quasi-Newton matrices}

\author{Jennifer B. Erway}
\email{erwayjb@wfu.edu}
\address{Department of Mathematics, PO Box 7388,
    Wake Forest University, Winston-Salem, NC 27109}
\author{Roummel F. Marcia}
\email{rmarcia@ucmerced.edu}
\address{School of Natural Sciences, University of California,
    Merced, 5200 N. Lake Road, Merced, CA 95343}
\thanks{Research supported in part by NSF grants CMMI-1334042 and CMMI-1333326.}

\begin{document}
\maketitle

\def\BFGS{{\small BFGS}}
\def\SR1{{\small SR1}}
\def\CG{{\small CG}}
\def\DFP{{\small DFP}}
\def\QR{{\small QR}}
\def\MATLAB{{\small MATLAB}}
\newcommand{\defined}{\mathop{\,{\scriptstyle\stackrel{\triangle}{=}}}\,}

\begin{abstract}
  In this paper, we consider the problem of efficiently computing the
  eigenvalues of limited-memory quasi-Newton matrices that exhibit a
  compact formulation.  In addition, we produce a compact formula for
  quasi-Newton matrices generated by any member of the Broyden convex class
  of updates.  Our proposed method makes use of efficient updates to the
  \QR{} factorization that substantially reduces the cost of
  computing the eigenvalues after the quasi-Newton matrix is updated.
  Numerical experiments suggest that the proposed method is able
    to compute eigenvalues to high accuracy.
  Applications for this work include modified quasi-Newton methods
  and trust-region methods for large-scale optimization, the efficient
  computation of condition numbers and singular values, and sensitivity
  analysis.
\end{abstract}


\pagestyle{myheadings}
\thispagestyle{plain}
\markboth{J. B. ERWAY AND R. F. MARCIA}{EIGENVALUES OF LIMITED-MEMORY QUASI-NEWTON MATRICES}

%%%%%%%%%%%%%%%%%%%%%%%%%%%%%%%%%%%%%%%%%%%%%%%%%%%%%%%%%%%%%%%%%%%%
\section{Introduction}

Newton's method for minimizing a twice-continuously differentiable
real-valued function $f:\Re^n\rightarrow \Re$ requires solving 
linear systems of the following form:
\begin{equation}\label{eq:Newton}
	\nabla^2 f(x_k) p_k = -\nabla f(x_k).
\end{equation}
When computing the Hessian is prohibitively expensive (e.g., $n$ is very
large) or $\nabla^2f(x_k)$ is unavailable, solving \eqref{eq:Newton} is
impractical or impossible.  In these cases, methods that only require
first-order information can be used in place of a pure Newton method.  Some
of the most commonly-used first-order methods are quasi-Newton methods,
which employ previous iterates and gradients to define a quasi-Newton
approximation $B_k$ to $\nabla^2 f(x_k)$, for each $k$.
% Using a matrix $B_k$ generated by
%quasi-Newton updates in place of the exact Hessian is common in both
%line-search and trust-region methods (citations).  
Conventional quasi-Newton matrices include Broyden-Fletcher-Goldfarb-Shanno
(\BFGS), Davidon-Fletcher-Powell (\DFP), symmetric rank-one (\SR1), and the
Broyden convex class.  In large-scale applications \emph{limited-memory}
versions of these quasi-Newton matrices are often used since they require
less storage.

In this paper we demonstrate how to efficiently compute the eigenvalues of
these limited-memory quasi-Newton matrices--allowing tools and methods
previously restricted to small or medium-scale problems to be used in a
large-scale quasi-Newton setting.  For example, the ability to efficiently
compute eigenvalues of these matrices allows for the use of traditional
modified Newton methods (see, e.g., \cite{GrNaS09, NocW99}).  In these methods, a \emph{modified} Newton step
is computed from solving a system of the form
\begin{equation*}
	\tilde{B}_k p_k = -\nabla f(x_k)
\end{equation*}
where $\tilde{B}_k$ is a positive-definite modification of $B_k$ obtained
by (i) replacing any negative eigenvalues with their absolute value in the spectral decomposition and (ii) thresholding to
%where $\tilde{B}_k$ is obtained from $B_k$ by modifying any negative
%eigenvalues so that $\tilde{B}_k$ is positive definite and thresholding to
prevent small eigenvalues. (For more details on modified Newton
  methods, please see, e.g., \cite[Section 11.4]{GrNaS09} or~\cite[Section
  3.4]{NocW99}.)  Knowledge of the eigenvalues of a quasi-Newton matrix
are also important in trust-region methods.  In these methods, one of the
challenges with using indefinite quasi-Newton matrices is in dealing with
the so-called ``hard case'', which occurs when the trust-region
  subproblem does not have a solution on the boundary of the trust region
  and $\nabla f(x_k)$ is perpendicular to the eigenspace associated with
  the most negative eigenvalue~\cite{ConGT00a,Gay81,MorS83}.  However,
knowing the leftmost eigenvalue provides an important bound in the ``hard
case'' that simplifies computations considerably (see,
e.g., %\cite{Johannes15} and
\cite[Algorithm 7.3.6]{ConGT00a}).

%both theoretical stability or convergence analysis or for devising .
Knowing the spectrum of the \BFGS, \DFP, \SR1{} or any member of the
Broyden convex class also enables the computation of their condition
numbers and singular values~\cite{GVL96}. From these, practical strategies
may be devised to avoid matrices that are poorly scaled or nearly singular.
Moreover, normwise sensitivity analysis for linear solves (e.g., forward
and backward error bounds) often use the condition number, the two-norm, or
the Frobenius norm of the system matrix (see e.g., \cite[Chapter
2.7]{GVL96} or \cite{Hig02}); with the results of this paper, these bounds
can be explicitly computed.

\medskip

% , computations of condition numbers can be performed for
% either practical computations adnd/or used for theoretical convergence
% analysis (see gvl?); Strategies can be devised that skip quasi-Newton
% updates can be used based on condition numbers to avoid matrices that are
% poorly scaled or nearly singular )kioday) Further, knowing the eigenvalues
% \emph{a priori} can lead to methods

%expands the scope of techniques and methods currently restricted for
%use wiht small-medieul size porblems to the large-scale arena.

%use for small to medium-scale problems applicable to large-scale optimization.

Our proposed method for efficiently computing eigenvalues of limited-memory
quasi-Newton matrices relies on compact representations of these matrices.
It is well-known that compact representations of \BFGS, \DFP, and \SR1
matrices are available \cite{BurWX96,ByrNS94,LDFP}; in particular, if $B$
is generated using any of these updates with an initial Hessian
approximation $B_0=\gamma I$, $\gamma\in\Re$, then $B$ can be written in
the form
\begin{equation}\label{eqn-compact0}
B=\gamma I + \Psi M \Psi^T,
\end{equation}
where $M$ is symmetric.  With this compact representation in hand, the
eigenvalue computation makes use of the \QR{} factorization of
$\Psi$.  This method was first proposed by Burdakov et
al.~\cite{Burdakov13}.  The bulk of the computational effort involves
computing the \QR{} factorization of $\Psi$.

This paper has two main contributions: (1) The compact representation of
the Broyden convex class of updates, and (2) the efficient updating of the
\QR{} factorization used for the eigenvalue decomposition.  We note that
while compact representations of \BFGS{} and \DFP{} matrices are
known~\cite{BurWX96,ByrNS94,LDFP}, to our knowledge there has been no work
done on generalizing the compact representation to the entire Broyden
convex class of updates.  This new compact representation allows us to
extend the eigenvalue computation to any member of the Broyden convex class
of updates.

\bigskip

Prior work on explicitly computing eigenvalues of quasi-Newton matrices has
been restricted to at most two updates.  A theorem by
Wilkinson~\cite[pp. 94--97]{Wil65} can be used to compute the eigenvalues
of a quasi-Newton matrix after one rank-one update.  For more than one
rank-one update, Wilkinson's theorem only provides bounds on
eigenvalues.  The eigenvalues of rank-one modifications are considered by
Golub~\cite{Golu73} and several methods are proposed, including Newton's
method on the characteristic equation, linear interpolation on a related
tridiagonal generalized eigenvalue problem, and finding zeros of the secular
equation.  Bunch et al.~\cite{BunNS78} extend the work of Golub to the
case of eigenvalue algebraic multiplicity of more than one. Eigenvalue
computations for more than one rank-one update are not proposed
and, as is, these methods cannot be used to compute the eigenvalues for the
general Broyden convex class of updates.
%. Lu and Monteiro~\cite{LuMo07} has 
%implementation when $B_k$ has special structure; namely, $B_k = D+VEV^T$,
%where $D$ and $E$ are positive diagonal matrices, and $V$ has a small
%number of columns. 
Apostolopoulou et al.~\cite{ApoSP08} and Apostolopoulou et
al.~\cite{ApoSP11} compute the eigenvalues of \emph{minimal}-memory \BFGS{}
matrices, where the number of \BFGS{} updates is limited to at most two.
In these papers, formulas for the characteristic polynomials are derived
that may be solved analytically.  Due to the complexity involved in
formulating characteristic polynomials and root finding, these approaches
cannot be generalized to handle more than two updates.  (In Appendix
  A, we show how the same characteristic polynomial for the case of one
  update as in~\cite{ApoSP11,ApoSP08} can be derived using our
  proposed approach.)

\bigskip

This paper is organized in six sections.  In Section 2, we outline the
compact formulations for the \BFGS, \DFP, and \SR1{} matrices.  In Section
3, we present the compact formulation for the Broyden convex class of
updates.  The method to compute the eigenvalues of any limited-memory
quasi-Newton matrix with the compact formulation (\ref{eqn-compact0}) is
given in Section 4.  An efficient method to update the \QR{} factorization
of $\Psi$ is also given in this section.  In Section 5 we demonstrate the
accuracy of the proposed method on a variety of limited-memory quasi-Newton
matrices.  Finally, in Section 6, there are some concluding remarks.

%%%%%%%%%%%%%%%%%%%%%%%%%%%%%%%%%%%%%%%%%%%%%%%%%%%%%%%%%%%%%%%%%%%%
\section{Compact formulations of quasi-Newton matrices}
In this section, we review compact formulations of some of the most
widely-used quasi-Newton matrices; in particular, we consider the \BFGS,
\DFP, and \SR1{} matrices.  First, we introduce
notation and assumptions used throughout this paper. 

\bigskip

Given a continuously differentiable function $f(x)\in\Re$ and
iterates $\{x_k\}$, the quasi-Newton pairs $\{s_i, y_i\}$ are defined as
follows: $$s_i \defined x_{i+1} - x_i \quad \text{and}\quad y_i \defined
\nabla f(x_{i+1}) - \nabla f(x_{i}),$$ where $\nabla f$ denotes the
gradient of $f$.  

\medskip

The goal of this section is to express a
quasi-Newton matrix obtained from these updates in the form
\begin{equation}\label{eqn-form}
	B_{k+1} \ = \ B_0 + \Psi_k M_k \Psi_k^T,
\end{equation}
where $\Psi_k \in \Re^{n \times l}$, $M_k \in \Re^{l \times l}$, and $B_0$
is a
diagonal matrix (i.e., $B_0=\gamma I$, $\gamma\in\Re$).  We will obtain
factorizations of the form (\ref{eqn-form}) where $l=k+1$ or $l=2(k+1)$; in
either case, we assume $l\ll n$.
%The form (\ref{eqn-form})
%has been used in various optimization settings \cite{BurWX96} (more?)  is
%already known for several (all?) of the updates considered in this paper
%(fix this depending on what we do).

\bigskip 

Throughout this section, we make use of the following matrices:
\begin{eqnarray*}
	S_k &\defined& [ \ s_0 \ \ s_1 \ \ s_2 \ \ \cdots \ \ s_{k} \ ] \ \in \ \Re^{n \times (k+1)}, \\
	Y_k &\defined& [ \ y_0 \ \ y_1 \ \ y_2 \ \ \cdots \ \ y_{k} \ ] \ \in \ \Re^{n \times (k+1)}.
\end{eqnarray*}
Furthermore, we make use of the following decomposition of $S_k^TY_k \in\Re^{(k+1) \times (k+1)}$:
$$
	S_k^TY_k =   L_k + D_k + R_k,
$$
where $L_k$ is strictly lower triangular, $D_k$ is diagonal, and $R_k$ is
strictly upper triangular.  We assume all updates are well-defined; for example,
for the \BFGS{} and \DFP{} updates, we assume that
$s_i^Ty_i > 0$ for $i = 0, 1, \dots, k$.
%Finally, for notational simplicity, we define
%$B_+\defined B_k$, $s \defined s_{k-1}$, and $y\defined y_{k-1}$.
\subsection{The BFGS update}  \label{sec-bfgs}
The Broyden-Fletcher-Goldfarb-Shanno (\BFGS) update 
is given by
\begin{eqnarray*}
	B_{k+1} &=& B_k  - \frac{1}{s_k^TB_k s_k}B_ks_k s_k^TB_k  +  
 	\frac{1}{y_k^Ts_k}y_ky_k ^T,
\end{eqnarray*}
where $B_0$ is a positive scalar multiple of the identity.  The \BFGS{}
update is the most widely-used rank-two update formula that (i) satisfies
the \emph{quasi-Newton condition} $B_{k+1}s_k=y_k$, (ii) has hereditary
symmetry, and (iii) provided that $y_i^Ts_i>0$ for $i=0,\ldots k$, then
$\{B_k\}$ exhibits hereditary positive-definiteness.  (For more background
on the \BFGS{} update formula, see, e.g., \cite{GrNaS09} or \cite{NocW99}.)

We now consider compact formulations of the \BFGS{} updates.
Byrd et al.~\cite[Theorem 2.3]{ByrNS94} showed that $B_{k+1}$ can be written in the form
\begin{equation}\label{eqn-alt-form}
	B_{k+1}= B_0 + \Psi_k\Gamma_k^{-1}\Psi_k^T,
\end{equation}
where
\begin{equation}\label{eqn-alt-form2}
	\Psi_k= \begin{pmatrix} B_0S_k & Y_k \end{pmatrix} \quad \text{and} \quad
	\Gamma_k = 
	-\begin{pmatrix}
	S_k^TB_0S_k & L_k \\
	L_k^T & -D_k
	\end{pmatrix}.
\end{equation}
Defining $M_k\defined \Gamma_k^{-1}$ gives us the desired compact form (\ref{eqn-form})
with $l=2(k+1)$.

\subsection{The DFP update}
The Davidon-Fletcher-Powell (\DFP) update is derived from applying \BFGS{}
updates to approximate the inverse of the Hessian.  The \DFP{} update
formula is given by
\begin{equation*}
	B_{k+1}  =
\left( I -\frac{y_ks_k^T}{y_k^Ts_k}\right) B_k
\left( I -\frac{s_ky_k^T}{y_k^Ts_k}\right) + \frac{y_ky_k^T}{y_k^Ts_k}
%	&=& B  - \frac{1 }{y ^Ts } y  s ^TB
%	- \frac{1}{y ^Ts }B s y ^T 
%	+ \frac{1}{(y ^Ts )^2}y  s ^TB  s  y ^T 
%	+ \frac{1}{y ^T s }y  y ^T.\\
\end{equation*}

%It can be shown that the \BFGS{} and \DFP{} formulas are related
%as follows:
%$$B_{k+1}^{DFP} = B_{k+1}^{BFGS} + (s_k^TB_ks_k)w_kw_k^T, \quad
%\text{where} \quad w = \frac{1}{y_k^Ts_k}y_k - \frac{1}{s_k^TB_ks_k}B_ks_k.$$ 
Like
\BFGS, \DFP{} updates satisfy the quasi-Newton condition while exhibiting
hereditary symmetry and positive-definiteness.  (For more background on the
\DFP{} update formula, see, e.g., \cite{GrNaS09} or \cite{NocW99}.)

The compact formulation for the \DFP{} update is found in
~\cite[Theorem 1]{LDFP}, where Erway et al. showed that 
$B_{k+1}$ can be written in the form (\ref{eqn-form})
with
\begin{equation}\label{eqn-dfpM}
	\Psi_k =
 	\begin{pmatrix}
       	B_0S_k & Y_k
     	\end{pmatrix}
	\quad \text{and} \quad
        M_k = 
	\begin{pmatrix}
		0 & -\bar{L}_k^{-T}\\
		- \bar{L}_k^{-1} & \bar{L}_k^{-1}(D_k+  S_k^TB_0S_k)\bar{L}_k^{-T}
	\end{pmatrix},
\end{equation}
where $\bar{L}_k\defined L_k+D_k$.  In this case, $l=2(k+1)$.
We note that in \cite{LDFP}, $\Psi_kM_k\Psi_k^T$ is expressed as
the equivalent product
\begin{equation*}
	\Psi_k M_k\Psi_k^T = 
	\begin{pmatrix}
       	Y_k & B_0S_k
     	\end{pmatrix}
	\begin{pmatrix}
		\bar{L}_k^{-1}(D_k+  S_k^TB_0S_k)\bar{L}_k^{-T} & -\bar{L}_k^{-1}\\
		- \bar{L}_k^{-T} & 0
	\end{pmatrix}
	\begin{pmatrix}
	Y_k^T\\
       	(B_0 S_k)^T 
     	\end{pmatrix}.
\end{equation*}
\subsection{The SR1 update}
The symmetric rank-one (\SR1) update 
formula is given by
\begin{eqnarray}\label{eqn-sr1}
	B_{k+1} &=& B_k  + \frac{1}
	{s_k^T(y_k  - B_k s_k ) }(y_k -B_k s_k )(y_k  - B_k s_k )^T.
\end{eqnarray}
The \SR1{} update is the unique rank-one update that satisfies the
quasi-Newton condition $B_{k+1}s_k=y_k$ and exhibits hereditary symmetry.
Unlike \BFGS{} and \DFP, these matrices do not exhibit hereditary
positive-definiteness. In fact, even if $y_i^Ts_i>0$ for each $i=0,\ldots
,k$, the sequence $\{B_k\}$ may not be positive definite.  The \SR1{}
update tends to be a better approximation of the true Hessian since it
is allowed to take on negative curvature; moreover, it has known 
convergence properties superior to other widely-used quasi-Newton
methods such as \BFGS{}~\cite{CGT91sr1}.
However, when using these updates extra precaution must be taken so the
denominator $s_k^T(y_k-B_ks_k)$ is nonzero.  (For more background on the
\SR1{} update formula, please see, e.g., \cite{NocW99} or \cite{GrNaS09}.)

The compact formulation for the \SR1{} update is due to 
Byrd et al.~\cite[Theorem 5.1]{ByrNS94}, where they showed that 
$B_{k+1}$
can be written in the form \eqref{eqn-form} with
%\begin{equation*}
%	 B_+ = \gamma I + 
%	\Psi	M^{-1} \Psi^T
%\end{equation*}
%where
\begin{equation*}
	\Psi_k \ = \ 
       	Y_k  - B_0S_k \quad \text{and} \quad
        M_k \ = \ (D_k + L_k + L_k^T - S_k^TB_0S_k)^{-1}.
\end{equation*}
Note that in the \SR1{} case, $l=k+1$.

\section{The Broyden convex class of updates}
In this section, we present a compact formulation for the Broyden
convex class of updates.
The Broyden convex class of updates is given by
 \begin{equation}\label{eqn-1param}
 B_{k+1}^\phi =B_k^\phi  -  \frac{1}{s_k ^TB_k^\phi s_k }B_k^\phi s_k s_k ^TB_k^\phi +  
  	\frac{1}{y_k^Ts_k }y_ky_k^T +
         \phi (s_k^T B_k^\phi s_k )w_kw_k^T,
 \end{equation}
 where $\phi \in [0, 1]$ and
$$
	w_k = \frac{y_k}{y_k^Ts_k} - \frac{B_k^\phi s_k}{s_k^TB_k^\phi s_k},
$$
(see, e.g., \cite{Luenberger,GrNaS09}).  Both the \BFGS{} and the \DFP{}
updates are members of this family.  (Setting $\phi=0$ gives the \BFGS{}
update, and setting $\phi=1$ yields the \DFP{} update.)  In fact, this
class of updates can be expressed in terms of the \BFGS{} and \DFP{}
updates: 
\begin{equation}\label{eqn-broyden}
  B_{k+1}^\phi=(1-\phi)B_{k+1}^{BFGS} + \phi B_{k+1}^{DFP},\end{equation}
where $\phi\in [0,1]$ and 
$$B_{k+1}^{BFGS}= 
B_k^\phi  - \frac{1}{s_k^TB_k^\phi s_k}B_k^\phi s_k s_k^TB_k^\phi  +  
 	\frac{1}{y_k^Ts_k}y_ky_k ^T,
$$
and
$$B_{k+1}^{DFP}=\left( I -\frac{y_ks_k^T}{y_k^Ts_k}\right) B_k^\phi
\left( I -\frac{s_ky_k^T}{y_k^Ts_k}\right) + \frac{y_ky_k^T}{y_k^Ts_k}.$$
(In other words, $B_{k+1}^\phi$ is the Broyden convex class matrix $B_k^\phi$ updated using
the \BFGS{} update and the \DFP{} update, respectively.)
All updates in this class satisfy the quasi-Newton condition.  Moreover,
members of this class enjoy hereditary symmetry and positive-definiteness 
provided $y_i^Ts_i>0$ for all $i$.  
Dixon~\cite{Dix72a} shows that, under some
conditions, the iterates generated using a quasi-Newton method belonging
to the Broyden class of convex updates together with an exact line search
will be identical in direction; in fact, the choice of $\phi$ only 
affects the step length of the search direction (see~\cite{Fle01} for
a detailed explanation).  In practice, for general nonlinear functions,
exact line searches are impractical; with inexact line searches,
it is well known that 
members of the Broyden convex class of updates can behave significantly
differently. (For more background and analysis of the Broyden convex
class of updates, see~\cite{Fle01}.)

We now consider compact formulations of the Broyden convex class of updates.
For notational simplicity, we drop the superscript $\phi$ for the duration
of this paper. 
To find the compact formulation, we expand (\ref{eqn-1param}) to obtain
\begin{eqnarray*}
	B_{k+1} &=& B_k - \frac{1-\phi}{s_k^TB_ks_k}B_ks_ks_k^TB_k 
			- \frac{\phi}{y_k^Ts_k}B_ks_ky_k^T 
			- \frac{\phi}{y_k^Ts_k}y_ks_k^TB_k \\
			&& \qquad + \left (1 + \phi\frac{s_k^TB_ks_k}{y_k^Ts_k} \right ) \frac{1}{y_k^Ts_k}y_ky_k^T,
\end{eqnarray*}
and thus, $B_{k+1}$ can be written compactly as
\begin{equation}\label{eqn-compact}
B_{k+1}=B_k + \left ( B_k s_k \ \ y_k \right )
		\begin{pmatrix}
			\displaystyle - \frac{(1-\phi)}{s_k^TB_ks_k} 
			&  \displaystyle - \frac{\phi}{y_k^Ts_k} \\
			\displaystyle -\frac{\phi}{y_k^Ts_k} 
			& \displaystyle  \left ( 1 + \phi \frac{s_k^TB_ks_k}{y_k^Ts_k}\right ) \frac{1}{y_k^Ts_k} 
		\end{pmatrix}
			\begin{pmatrix}
				(B_ks_k)^T \\
				y_k^T
			\end{pmatrix}.
\end{equation}
Recall that our goal is to write $B_{k+1}$ in the form
$$B_{k+1}=B_0+\Psi_k M_k\Psi_k^T,$$
where $
\Psi_k 
\in  \Re^{n \times 2(k+1)}$ and $M_k \in \Re^{2(k+1) \times 2(k+1)}$. 
Letting $\Psi_k$ be defined as 
\begin{equation}\label{eqn-psik} 
%\Psi_k \defined [ \ B_0s_0 \ \ B_0s_1 \ \ \cdots \ \ B_0s_{k} \ \ y_0 \ \ y_1 \ \ \cdots \ \ y_k \ ] = \begin{pmatrix} \ B_0S_k & Y_k \end{pmatrix}
\Psi_k \defined \begin{pmatrix} B_0s_0 & B_0s_1 & \cdots & B_0s_{k} & y_0 & y_1 & \cdots & y_k \end{pmatrix} = \begin{pmatrix} \ B_0S_k & Y_k \end{pmatrix},
\end{equation}
we now derive an expression for $M_k$.

\bigskip

\subsection{General $M_k$} \label{eq:GeneralG} In this section we state and
prove a theorem that gives an expression for $M_k$.  The eigenvalue
computation in Section 5 requires the ability to form $M_k$.  For this
reason, we also provide a practical recursive method for computing $M_k$.  \bigskip

\noindent 
\textbf{Theorem 1.} 
\textsl{
%Let $S_k = [s_0 \ \ s_1 \ \ \cdots \ \ s_k] \in \Re^{n
%  \times (k+1)}$ and $Y_k = [y_0 \ \ y_1 \ \ \cdots \ \ y_k] \in \Re^{n
%  \times (k+1)}$, and let $S_k^TY_k= L_k + D_k + R_k$, where $L_k \in
%\Re^{(k+1) \times (k+1)}$ is strictly lower triangular, $D_k \in \Re^{(k+1)
%  \times (k+1)}$ is diagonal and $R_k \in \Re^{(k+1) \times (k+1)}$ is
%strictly upper triangular.  
Let $\Lambda_k \in \Re^{(k+1) \times (k+1)}$ be
a diagonal matrix such that
\begin{equation}\label{eq:Lambda}
	\Lambda_k =  \underset{0 \le i \le k}{\text{diag}} \big ( \lambda_i \big ),
	\qquad \text{where \ } 
			\lambda_i = 
			\frac{1}{\displaystyle 
			-\frac{1-\phi}{s_i^TB_is_i}
			-\frac{\phi}{s_i^Ty_i}}
	\ \text{for $0 \le i \le k$.}
\end{equation}
If $B_{k+1}$ is updated using the Broyden convex class of updates (\ref{eqn-1param}), where $\phi \in [0,1]$, then $B_{k+1}$ can
be written as $B_{k+1}=B_0+\Psi_k M_k\Psi_k^T$, where $\Psi_k$ is defined
as in (\ref{eqn-psik}) and
\begin{equation}\label{eqn-G} 
M_k=	\begin{pmatrix}
		-S_k^TB_0S_k + \phi \Lambda_k  & -L_k + \phi \Lambda_k \\
		-L_k^T + \phi \Lambda_k & \ \ D_k + \phi \Lambda_k
	\end{pmatrix}^{-1}.
\end{equation}
}

\bigskip

\noindent \textbf{Proof.} This proof is broken into two parts.
First, we consider the special cases when $\phi=0$ and $\phi=1$.  Then, we
prove by induction the case when $\phi\in (0,1)$.

When $\phi = 0$,  (\ref{eqn-1param}) becomes the \BFGS{} update
and $M_k$ in \eqref{eqn-G} simplifies to
$$
	M_k = 
	\begin{pmatrix}
		-S_k^TB_0S_k & -L_k \\
		-L_k^T & \ \ \ D_k
	\end{pmatrix}^{-1},
$$ 
which is consistent with (\ref{eqn-alt-form}) and (\ref{eqn-alt-form2}).
When $\phi = 1$, then $\Lambda_k = -D_k$ 
and so
\eqref{eqn-1param} is the \DFP{} update and
with
$$
	M_k = 
	\begin{pmatrix}
		-S_k^TB_0S_k-D_k & -\bar{L}_k \\
		-\bar{L}_k^T & 0
	\end{pmatrix}^{-1},
$$ 
where $\bar{L} = L_k + D_k$.  After some algebra, 
it can be shown that this is exactly $M_k$ given in (\ref{eqn-dfpM}).
Thus, $M_k$ in \eqref{eqn-G} is correct for
$\phi = 0$ and $\phi = 1$.

\medskip

The proof for $\phi\in (0,1)$ is by induction on $k$.
We begin by considering the base case $k=0$.
For $k=0$, $B_1$ is given by (\ref{eqn-compact}), and thus,
$B_1=B_0+\Psi_0 \widehat{M}_0\Psi_0^T$ where $\Psi_0 = ( \ B_0s_0 \ \ y_0 \ )$
and  \begin{equation}\label{eq:G0}
	\widehat{M}_0 \defined
		\begin{pmatrix}
			\displaystyle - \frac{(1-\phi)}{s_0^TB_0s_0} 
			&  \displaystyle - \frac{\phi}{y_0^Ts_0} \\
			\displaystyle -\frac{\phi}{y_0^Ts_0} 
			& \displaystyle  \left ( 1 + \phi \frac{s_0^TB_0s_0}{y_0^Ts_0}\right ) \frac{1}{y_0^Ts_0} 
		\end{pmatrix}.
\end{equation}
To complete the base case, we now show that $\widehat{M}_0$ in (\ref{eq:G0}) is
equivalent to $M_0$ in (\ref{eqn-G}).
For simplicity, $\widehat{M}_0$ can be written as
\begin{equation}	
\widehat{M}_0	=
		\begin{pmatrix}
			\alpha_0 & \beta_0 \\
			\beta_0 & \delta_0
		\end{pmatrix},
\end{equation}
where
\begin{equation}\label{eqn-alphabetadelta}
	\alpha_0 = - \frac{(1-\phi)}{s_0^TB_0s_0}, \quad
	\beta_0 = - \frac{\phi}{y_0^Ts_0}, \quad \text{and} \quad 
	\delta_0  = \left ( 1 + \phi \frac{s_0^TB_0s_0}{y_0^Ts_0}\right ) \frac{1}{y_0^Ts_0}.
\end{equation}
We note that $\alpha_0$ and $\beta_0$ are nonzero  since $0 < \phi < 1$.  
Consequently, $\delta_0$ can be written as
%\begin{eqnarray}
%\delta_0 &=& \left ( 1 + \phi \frac{s_0^TB_0s_0}{y_0^Ts_0}\right )
%\frac{1}{y_0^Ts_0} \nonumber
%\\ &= &
%	- \left ( 1 + (1-\phi) \frac{\beta_0}{\alpha_0} \right ) \frac{\beta_0}{\phi} \nonumber
%\\  & = &
%	-\frac{\beta_0}{\phi} - \frac{\beta_0^2}{\phi\alpha_0} + \frac{\beta_0^2}{\alpha_0}. \label{eqn-delta}
%\end{eqnarray}
\begin{equation}
\delta_0 = \left ( 1 + \phi \frac{s_0^TB_0s_0}{y_0^Ts_0}\right )
\frac{1}{y_0^Ts_0} =
	- \left ( 1 + (1-\phi) \frac{\beta_0}{\alpha_0} \right ) \frac{\beta_0}{\phi} =
	-\frac{\beta_0}{\phi} - \frac{\beta_0^2}{\phi\alpha_0} + \frac{\beta_0^2}{\alpha_0}. \label{eqn-delta}
\end{equation}
The determinant, $\eta_0$, of $\widehat{M}_0$ can be written as 
\begin{eqnarray}\label{eqn-eta0}
\eta_0=	\alpha_0 \delta_0 - \beta_0^2 \ = \
	-\frac{\alpha_0 \beta_0}{\phi} - \frac{\beta_0^2}{\phi} \ = \
	-\frac{\beta_0}{\phi}(\alpha_0 + \beta_0).
\end{eqnarray}
Since all members of the convex class are positive 
definite, %$B_0 = \gamma I$ is positive definite and
%$s_0^Ty_0>0$, 
both $\alpha_0$ and $\beta_0$ are negative,
and thus, $\alpha_0+\beta_0\ne 0$ and $\eta_0 \ne 0$ in \eqref{eqn-eta0}.
It follows that $\widehat{M}_0$ is invertible, and in particular,
$$
	\widehat{M}_0^{-1} = 
	\begin{pmatrix}
		\delta_0/\eta_0 & - \beta_0/\eta_0 \\
		-\beta_0/\eta_0 & \alpha_0/\eta_0
	\end{pmatrix}.
$$
%where $\eta_0\defined \alpha_0\delta_0 - \beta_0^2$. %\frac{1}{\alpha_0\delta_0 - \beta_0^2}$.
%Using (\ref{eqn-delta}), $\eta_0$ can be simplified as follows:
Together with (\ref{eqn-delta}), 
the (1,1) entry of $\widehat{M}_0^{-1}$ simplifies to
\begin{eqnarray}
  \frac{\delta_0}{\eta_0} &=&
  \frac{\displaystyle
    -\left (
      \frac{\alpha_0 + \beta_0}{\alpha_0} - \frac{ \phi \beta_0}{\alpha_0} 
    \right )
    \frac{\beta_0}{\phi}
  }{\displaystyle - \frac{\beta_0}{\phi}(\alpha_0 + \beta_0)} \nonumber \\ 
  &=&
  \frac{1}{\alpha_0} -\frac{\phi\beta_0}{\alpha_0 (\alpha_0 + \beta_0)}\nonumber \\[.3cm]
%\frac{1}{\alpha_0 + \beta_0} \\[.3cm] 
	&=&
	\frac{(\alpha_0+\beta_0)(1 - \phi)+\phi \alpha_0}{\alpha_0(\alpha_0+\beta_0)}\nonumber \\[.3cm]
	&=&
	\frac{1-\phi}{\alpha_0}  + \frac{\phi}{\alpha_0+\beta_0}. \label{eqn-deltaeta}
\end{eqnarray}
Finally, the (2,2) entry of $\widehat{M}_0^{-1}$ can be written as
\begin{eqnarray}\label{eqn-22}
	\frac{\alpha_0}{\eta_0} &=& -\frac{\phi \alpha_0}
	{\beta_0(\alpha_0+\beta_0)} \ = \ 
	-\frac{\phi}{\beta_0} + \frac{\phi}{\alpha_0+\beta_0}.
\end{eqnarray}
Thus, combining (\ref{eqn-eta0}), (\ref{eqn-deltaeta}),  and (\ref{eqn-22}),
we obtain the following equivalent expression for $\widehat{M}_0^{-1}$:
\begin{eqnarray}\label{eqn-M0inv}
	\widehat{M}_0^{-1} &=&
	\begin{pmatrix}
		\displaystyle  \frac{1-\phi}{\alpha_0}  + \frac{\phi}{\alpha_0+\beta_0} &
		\displaystyle   \frac{\phi}{\alpha_0 + \beta_0}	
		\\[.4cm]
		\displaystyle  \frac{\phi}{\alpha_0 + \beta_0} &
		\displaystyle  -\frac{\phi}{\beta_0} + \frac{\phi}{\alpha_0+\beta_0}
	\end{pmatrix}.
\end{eqnarray}
% 	\ = \
% 	\begin{pmatrix}
% 		s_0^TB_0s_0 + \phi \lambda_0 & \phi \lambda_0 \\
% 		\phi \lambda_0 & s_0^Ty_0 + \phi \lambda_0
% 	\end{pmatrix},
% \end{eqnarray*}
For the case $k=0$, $L_0=R_0=0$; moreover, 
$\lambda_0 = 1/(\alpha_0 + \beta_0)$.
Substituting back in for $\alpha_0$, $\beta_0$ and $\delta_0$ using
(\ref{eqn-alphabetadelta}), we
obtain
\begin{eqnarray*}
	\widehat{M}_0  &=&
	\begin{pmatrix}
		-s_0^TB_0s_0 + \phi \lambda_0 & \phi \lambda_0 \\
		\phi \lambda_0 & s_0^Ty_0 + \phi \lambda_0
	\end{pmatrix}^{-1},
\\
& = &
	\begin{pmatrix}
		-S_k^TB_0S_k + \phi \Lambda_k  & -L_k + \phi \Lambda_k \\
		-L_k^T + \phi \Lambda_k & \ \ D_k + \phi \Lambda_k
	\end{pmatrix}^{-1} \\
& = & 
M_0,
\end{eqnarray*}
proving the base case.

\bigskip

For the induction step, assume 
\begin{equation}\label{eqn-inductB}
B_{m}=B_0+\Psi_{m-1} M_{m-1}\Psi_{m-1}^T,
\end{equation} where $\Psi_{m-1}$ is defined
as in (\ref{eqn-psik}) and 
\begin{eqnarray}\label{eqn-Mm-1}
M_{m-1}=	\begin{pmatrix}
		-S_{m-1}^TB_0S_{m-1} + \phi \Lambda_{m-1}  & -L_{m-1} + \phi \Lambda_{m-1} \\
		-L_{m-1}^T + \phi \Lambda_{m-1} & \ \ D_{m-1} + \phi \Lambda_{m-1}
	\end{pmatrix}^{-1}.
\end{eqnarray}
From (\ref{eqn-compact}), we have
\begin{eqnarray}\label{eqn-bnp1}
        \ \ B_{m+1} &=& B_0 + \Psi_{m-1}M_{m-1}\Psi_{m-1}^T +
                        \left ( B_ms_m \ \ \ y_m \right )
	                \begin{pmatrix}
                                \alpha_m & \beta_m \\
	                        \beta_m & \delta_m
                        \end{pmatrix}
             	        \begin{pmatrix}
                                (B_ms_m)^T \\
   	                        y_m^T
                        \end{pmatrix},
\end{eqnarray}
where
\begin{eqnarray*}
        \alpha_m = -\frac{1-\phi}{s_m^TB_ms_m}, \quad
        \beta_m  = - \frac{\phi}{y_m^Ts_m},
\end{eqnarray*}
and
\begin{eqnarray*}
        \delta_m  =   \left ( 1 + \phi \frac{s_m^TB_ms_m}{y_m^Ts_m}\right ) \frac{1}{y_m^Ts_m}  =
         \displaystyle - \left ( 1 + (1-\phi)\frac{\beta_m}{\alpha_m} \right ) \
\frac{\beta_m}{\phi}.
\end{eqnarray*}
As in the base case, $k = 0$, we note that $\alpha_m$ and $\beta_m$ are nonzero since $0 < \phi < 1$,
and that the determinant $\alpha_m \delta_m - \beta_m^2$ is also nonzero.

Multiplying (\ref{eqn-inductB}) by $s_m$ on the right, we obtain
\begin{equation}\label{eqn-block22}
B_ms_m = B_0s_m + \Psi_{m-1}M_{m-1}\Psi_{m-1}^Ts_m.
\end{equation}
Then, substituting this into (\ref{eqn-bnp1}) yields
\begin{equation}\label{eqn-bnp12}
	B_{m+1}
	= B_0 + \Psi_{m-1}M_{m-1}\Psi_{m-1}^T +
			\left ( B_0s_m + \Psi_{m-1}p_m \ \ \ y_m \right )
			\begin{pmatrix}
				\alpha_m & \beta_m \\
				\beta_m & \delta_m
			\end{pmatrix}
			\begin{pmatrix}
				(B_0s_m + \Psi_{m-1}p_m)^T \\
				y_m^T
			\end{pmatrix},
\end{equation}
where $p_m\defined M_{m-1}\Psi_{m-1}^Ts_m$.
Equivalently,
\begin{equation} \label{eqn-3x3}
	B_{m+1}
	= B_0 +  \left ( \Psi_{m-1} \ \ B_0s_m \ \ y_m \right )
		\begin{pmatrix}
			M_{m-1}+ \alpha_m p_mp_m^T & \alpha_m p_m & \beta_m p_m \\
			\alpha_m p_m^T & \alpha_m & \beta_m \\
			\beta_m  p_m^T & \beta_m & \delta_m
		\end{pmatrix}
		\begin{pmatrix}
			\Psi_{m-1}^T \\
			(B_0s_m)^T \\
			y_m^T
		\end{pmatrix}.
\end{equation}

The $3\times 3$ block matrix in (\ref{eqn-3x3}) has the following decomposition:
\begin{equation}\label{eq:GkP}
		\begin{pmatrix}
			M_{m-1}+ \alpha_m p_mp_m^T & \alpha_m p_m & \beta_m p_m\\
			\alpha_m  p_m ^T & \alpha_m  & \beta_m \\
			\beta_m  p_m ^T & \beta_m  & \delta_m 
		\end{pmatrix}
		=
		\begin{pmatrix}
			I & p_m  & 0 \\
			0 & 1 & 0 \\
			0 & 0 & 1
		\end{pmatrix}
		\begin{pmatrix}
			M_{m-1} & 0 & 0 \\
			0 & \alpha_m  & \beta_m  \\
			0 & \beta_m  & \delta_m 
		\end{pmatrix}
		\begin{pmatrix}
			I & 0 & 0 \\
			p_m ^T & 1 & 0 \\
			0 & 0 & 1
		\end{pmatrix},
\end{equation}
allowing us to compute its inverse as follows:
\begin{eqnarray}\nonumber
		\begin{pmatrix}
			M_{m-1}+ \alpha_m p_mp_m^T & \alpha_m p_m & \beta_m p_m\\
			\alpha_m  p_m ^T & \alpha_m  & \beta_m \\
			\beta_m  p_m ^T & \beta_m  & \delta_m 
		\end{pmatrix}^{-1}
		&=&
		\begin{pmatrix}
			I & 0 & 0 \\
			-p_m^T & 1 & 0 \\
			0 & 0 & 1
		\end{pmatrix}
		\begin{pmatrix}
			M_{m-1}^{-1} & 0 & 0 \\
			0 & \tilde{\alpha}_m & \tilde{\beta}_m \\
			0 & \tilde{\beta}_m &  \tilde{\delta}_m
		\end{pmatrix}
		\begin{pmatrix}
			I & -p_m & 0 \\
			0 & 1 & 0 \\
			0 & 0 & 1
		\end{pmatrix}
		\\
		&=&
		\begin{pmatrix}
		M_{m-1}^{-1} & -M_{m-1}^{-1}p_m & 0 \\
		-p_m^TM_{m-1}^{-1} & p_m^TM_{m-1}^{-1}p_m  + \tilde{\alpha}_m & \tilde{\beta}_m \\
		0 & \tilde{\beta}_m & \tilde{\delta}_m
		\end{pmatrix}, \label{eq:Ginv}
\end{eqnarray}
where
$$
\begin{array}{rlllllll}
  \tilde{\alpha}_m &=& \displaystyle \frac{\delta_m}{\alpha_m \delta_m - \beta_m^2} &=&\displaystyle 
\frac{1-\phi}{\alpha_m}  + \frac{\phi}{\alpha_m+\beta_m}\\
  \tilde{\beta}_m &=& \displaystyle \frac{-\beta_m}{\alpha_m \delta_m - \beta_m^2} &=&
  \displaystyle \frac{\phi}{\alpha_m+\beta_m} \\
\displaystyle  \tilde{\delta}_m &=& \displaystyle \frac{\alpha_m}{\alpha_m \delta_m - \beta_m^2} &=&
 \displaystyle 	-\frac{\phi}{\beta_0} + \frac{\phi}{\alpha_0+\beta_0}.
\end{array}
$$

We now simplify the entries of (\ref{eq:Ginv}).
Since $p_m = M_{m-1}\Psi_{m-1}^Ts_m$, then $M_{m-1}^{-1}p_m = \Psi_{m-1}^Ts_m$, giving us an expression for the (1,2) and (2,1) entries.  The (2,2) block entry 
is simplified by first multiplying (\ref{eqn-block22}) by $s_m^T$ on the left to obtain $s_m^TB_ms_m = s_m^TB_0s_m + p_m^TM_{m-1}^{-1}p_m$.  Then,
\begin{eqnarray*}
	p_m^TM_{m-1}^{-1}p_m + \tilde{\alpha}_m &=&  -s_m^TB_0s_m + s_m^TB_ms_m + \tilde{\alpha}_m \\
	&=& -s_m^TB_0s_m - \frac{1-\phi}{\alpha_m} +
	\frac{1-\phi}{\alpha_m}  + \frac{\phi}{\alpha_m+\beta_m}\\
	&=& -s_m^TB_0s_m + \frac{\phi}{\alpha_m + \beta_m}.
\end{eqnarray*}
Thus, (\ref{eq:Ginv}) can be written as
\begin{equation}\label{eqn-Ginv-simp}
		\begin{pmatrix}
		M_{m-1}^{-1} & -\Psi_{m-1}^Ts_m & 0 \\
		-s_m^T\Psi_{m-1} & -s_m^TB_0s_m\! +\! \displaystyle \frac{\phi}{\alpha_m+\beta_m}  
		& \displaystyle \frac{\phi}{\alpha_m + \beta_m} \\
		0 & \displaystyle \frac{\phi}{\alpha_m + \beta_m} 
		& y_m^Ts_m\! +\! \displaystyle \frac{\phi}{\alpha_m+\beta_m}
		\end{pmatrix}.	
\end{equation}

\medskip

We now show (\ref{eqn-G}) holds using (\ref{eqn-Ginv-simp}).
%psi n is n x 2(n+1)
Define the permutation matrix $\Pi_m\in\Re^{2(m+1)\times 2(m+1)}$ as follows:
\begin{equation}\label{eqn-permute}
	\Pi_m = 
	\begin{pmatrix}
		I_m  \\
		& 0& I_m \\
                & 1 & 0 \\
		& & & 1
	\end{pmatrix}.
\end{equation}
Then,
$$	\left [ \Psi_{m-1} \ \ B_0s_m \ \ y_m \right ]
	\Pi_m
	= \Psi_m;
$$
in other words, $[ \Psi_{m-1} \ \ B_0s_m \ \ y_m] = \Psi_m \Pi_m^T$.  Therefore, (\ref{eqn-3x3}) can be written as
$$
	B_{m+1}
	= B_0 +
	\Psi_m\Pi_m^T\widehat{M}_m	\Pi_m \Psi_m^T, 
$$
where
$$\widehat{M}_m\defined
\begin{pmatrix}
			M_{m-1}+ \alpha_m p_mp_m^T & \alpha_m p_m & \beta_m p_m \\
			\alpha_m p_m^T & \alpha_m & \beta_m \\
			\beta_m p_m^T & \beta_m & \delta_m
		\end{pmatrix}.
$$
It remains to show (\ref{eqn-G}) holds for $k=m$ with
$M_m\defined \Pi_m^T\widehat{M}_m\Pi_m$.

\medskip
We now consider $M_m^{-1}$ given by
\begin{eqnarray}\label{eqn-Mminv}
	M_m^{-1} &=& 
	\left (
	\Pi_m^T 
		\begin{pmatrix}
			M_{m-1}+ \alpha_m p_mp_m^T & \alpha_m p_m & \beta_m p_m \\
			\alpha_m p_m^T & \alpha_m & \beta_m \\
			\beta_m p_m^T & \beta_m & \delta_m
		\end{pmatrix}
	\Pi_m
	\right )^{-1},
\end{eqnarray}
which can be simplified using 
(\ref{eqn-Ginv-simp}):
\begin{equation}\label{eqn-Minv}
	M_m^{-1}= \Pi_m^T
		\begin{pmatrix}
		M_{m-1}^{-1} & -\Psi_{m-1}^Ts_m & 0 \\
		-s_m^T\Psi_{m-1} & -s_m^TB_0s_m\! +\! \displaystyle \frac{\phi}{\alpha_m+\beta_m}  
		& \displaystyle \frac{\phi}{\alpha_m + \beta_m} \\
		0 & \displaystyle \frac{\phi}{\alpha_m + \beta_m} 
		& y_m^Ts_m\! +\! \displaystyle \frac{\phi}{\alpha_m+\beta_m}
		\end{pmatrix}
	\Pi_m.
\end{equation}

Now partition $M_{m-1}^{-1}$ as follows:
$$
	M_{m-1}^{-1} \ = \
	\begin{pmatrix}
		(M_{m-1}^{-1})_{11} & (M_{m-1}^{-1})_{12} \\
		(M_{m-1}^{-1})_{21} & (M_{m-1}^{-1})_{22}
	\end{pmatrix}.
$$	
Applying the permutation matrices together with $\Psi_{m-1}^Ts_m = \displaystyle
\begin{pmatrix}
	S_{m-1}^TB_0s_m \\
	Y_{m-1}^Ts_m
\end{pmatrix}
$, we have that
\begin{equation*}
M_m^{-1}=
	\left (
	\begin{tabular}{cc|ccccc}
		$(M_{m-1}^{-1})_{11}$ & $-S_{m-1}^TB_0s_m$ & $(M_{m-1}^{-1})_{12}$ & $0$ \\
		$-s_m^TB_0S_{m-1}$& $-s_m^TB_0s_m+ \displaystyle \frac{\phi}{\alpha_m+\beta_m}$ &
		$-s_m^TY_m$ & $\displaystyle \frac{\phi}{\alpha_m + \beta_m}$  \\\hline
		$(M_{m-1}^{-1})_{21}$ & $-Y_m^Ts_m$ & $(M_{m-1}^{-1})_{22}$ & $0$ \\
		$0$ & $\displaystyle \frac{\phi}{\alpha_m + \beta_m}$ & $0$ & $y_m^Ts_m\! +\! \displaystyle \frac{\phi}{\alpha_m+\beta_m}$
	\end{tabular}
	\right ).
\end{equation*}
Simplifying using the induction hypothesis (\ref{eqn-Mm-1}) yields
\begin{eqnarray*}
M_m^{-1} & = &
	\left (
	\begin{tabular}{cc|ccccc}
		$-S_{m-1}^TB_0S_{m-1} + \phi \Lambda_{m-1}$ & $-S_{m-1}^TB_0s_m$ & $-L_{m-1} + \phi \Lambda_{m-1}$ & $0$ \\
		$-s_m^TB_0S_{m-1}$& $-s_m^TB_0s_m+ \displaystyle \frac{\phi}{\alpha_m+\beta_m}$ &
		$-s_m^TY_m$ & $\displaystyle \frac{\phi}{\alpha_m + \beta_m}$  \\[.2cm] \hline
		$-L_{k-1}^T+\phi \Lambda_{k-1}$ & $-Y_m^Ts_m$ & $D_{k-1} + \phi \Lambda_{k-1}$ 
		& $0$ \\
		$0$ & $\displaystyle \frac{\phi}{\alpha_m + \beta_m}$ & $0$ & $y_m^Ts_m\! +\! \displaystyle \frac{\phi}{\alpha_m+\beta_m}$
	\end{tabular}
	\right )\\
	&=& 	
	\begin{pmatrix}
		-S_m^TB_0S_m + \phi \Lambda_m  & -L_m + \phi \Lambda_m \\
		-L_m^T + \phi \Lambda_m & \ \ D_m + \phi \Lambda_m
	\end{pmatrix},
\end{eqnarray*}
i.e., (\ref{eqn-G}) holds for $k=m$.  $\square$

\bigskip

Although we have found an expression for $M_k$, computing $M_k$ is not
straightforward.  In particular, the diagonal matrix $\Lambda_k$ in Eq.\
\eqref{eq:Lambda} involves $s_i^TB_is_i$, which requires $B_i$ for $0 \le i
\le k$.  In the following section we propose a different way of computing
$M_k$ that does not necessitate storing the quasi-Newton matrices $B_i$ for
$0 \le i \le k$.

\bigskip

\subsection{Computing $M_k$} \label{sec:ComputingM} In this section we
propose a recursive method for computing $M_k$ from $M_{k-1}$.
We have shown that
$$
	M_k \ = \ 
	\Pi_k^T
	\begin{pmatrix}
			M_{k-1}+ \alpha_k p_kp_k^T & \alpha_k p_k & \beta_k p_k \\
			\alpha_k p_k^T & \alpha_k & \beta_k \\
			\beta_k  p_k^T & \beta_k & \delta_k
	\end{pmatrix}
	\Pi_k.
$$

The vector $p_k$ can be computed
as follows:
$$
	p_k = M_{k-1}\Psi_{k-1}^Ts_k =
	M_{k-1}
	\begin{pmatrix}
		(B_0S_{k-1})^T \\
		Y_{k-1}^T
	\end{pmatrix}
	s_k
	=
	M_{k-1}
	\begin{pmatrix}
		S_{k-1}^TB_0s_k \\
		Y_{k-1}^Ts_k
	\end{pmatrix}.
$$
Note that $(S_{k-1}^TB_0s_k)^T$ is the last row (save the diagonal entry)
of $S_k^TB_0S_k$ and $(Y_{k-1}^Ts_k)^T$ is the last row (save the diagonal
entry) of $S_k^TY_k$.  The entry $\alpha_k$, which
is given by $\alpha_k = -(1-\phi)/s_k^TB_ks_k$ can be
computed from the following:
\begin{equation}\label{eqn-makingalpha}
	s_k^TB_ks_k = s_k^T \bigg (B_0 + \Psi_{k-1}M_{k-1}\Psi_{k-1}^T \bigg ) s_k
	= s_k^TB_0s_k + s_k^T\Psi_{k-1}p_k.
\end{equation}
The quantity $s_k^TB_0s_k$ is the last diagonal entry in $S_k^TB_0S_k$, and
$s_k^T\Psi_{k-1}p_k$ is the inner product of $\Psi_{k-1}^Ts_k$ (which was
formed when computing $p_k$) and $p_k$.  The entry $\beta_k$ is given by
$\beta_k = -\phi/y_k^Ts_k$, where $y_k^Ts_k$ is the last diagonal entry in
$S_k^TY_k$.  Finally, $\delta_k = ( 1 + \phi s_k^TB_ks_k/y_k^Ts_k) / y_k^Ts_k
$, which uses the previously computed quantities $s_k^TB_ks_k$
and $y_k^Ts_k$.

\medskip

We summarize this recursive method in Algorithm 1.

\bigskip

\noindent \textbf{Algorithm 1.} This algorithm computes $M_k$ in (\ref{eqn-G}).

\bigskip

Define $\phi$ and $B_0$;

Define $M_0$ using (\ref{eq:G0});

Define $\Psi_0 = (B_0s_0 \,\, y_0 )$;

\textbf{for} $j = 1:k$

\quad $p_j \gets M_{j-1} (\Psi_{j-1}^Ts_j$);

\quad $s_j^TB_js_j \gets s_j^TB_0s_j + ( s_j^T\Psi_{j-1})p_j$;

\quad $\alpha_j \gets -(1-\phi)/ (s_j^TB_js_j$);

\quad $\beta_j \gets -\phi/ (y_j^Ts_j$);

\quad $\delta_j \gets  ( 1 + \phi (s_j^TB_js_j)/(y_j^Ts_j) ) / (y_j^Ts_j)$;

\quad 
$\displaystyle 
	M_j \ \gets \ 
	\Pi_j^T
	\begin{pmatrix}
			M_{j-1}+ \alpha_j p_jp_j^T & \alpha_j p_j & \beta_j p_j  \\
			\alpha_j p_j^T & \alpha_j & \beta_j \\
			\beta_j  p_j^T & \beta_j & \delta_j
	\end{pmatrix}
	\Pi_j
$, where $\Pi_j$ is as in (\ref{eqn-permute});

\textbf{end}

\bigskip 

The matrices $\{\Pi_j\}$ are not formed explicitly since they
are permutation matrices; thus, no matrix-matrix products are required
by the recursion algorithm.

%%%%%%%%%%%%%%%%%%%%%%%%%%%%%%%%%%%%%%%%%%%%%%%%%%%%%%%%%%%%%%%%%%%%
\section{Computing the eigenvalues of $B_{k+1}$}
In this section, we demonstrate how to compute the eigenvalues of
a limited-memory matrix
$B_{k+1}$ when the following decomposition is available:
\begin{equation}\label{eqn-compactDe}
	B_{k+1} = B_0 + \Psi_k M_k \Psi_k^T,
\end{equation}
where $B_0 = \gamma I$, $\gamma\in \Re$.  We assume that
$B_{k+1}\in\Re^{n\times n}$ but only $m$ limited-memory updates are stored,
where $m\ll n$ (see, e.g., \cite{kolda1998bfgs, liu1989limited,
  nash1991numerical, nocedal1980updating}).  In large-scale optimization,
typically $m<10$ (e.g., Byrd et al.~\cite{ByrNS94} recommend $m\in [2,6]$).

\medskip

The material presented in Section 4.1 was first proposed by Burdakov et
al.~\cite{Burdakov13} in a different manner.  We explain
these differences at the end of the section.

\subsection{Eigenvalues via the QR decomposition}
We begin by finding the eigenvalues of $B_{k+1}$ when $B_{k+1}$ is obtained
using the Broyden convex class of updates; at the end of this section, we
describe the modifications needed to find the eigenvalues for the \SR1{}
case.  We assume $k+1\le m \ll n$.

For the Broyden convex class of updates, $\Psi_k=\begin{pmatrix} B_0 S_k &
  Y_k \end{pmatrix}$, i.e.,
$$\Psi_k =
\begin{pmatrix}
B_0 s_0 & B_0s_1 & \cdots & B_0s_{k} & y_0 & y_1 & \cdots & y_{k}
\end{pmatrix}.
$$
To facilitate updating $\Psi_k$ after computing a new limited-memory pair
(see Section~\ref{sec-updates}), we permute the columns of $\Psi_k$ using a
permutation matrix $P$ (also called the ``perfect shuffle''--see,
e.g., \cite{kressner}) so that
$$\hat{\Psi}_k \defined \Psi_k P =
\begin{pmatrix}
B_0 s_0 & y_0 & B_0 s_1 & y_1 & \cdots & B_0 s_{k} & y_{k}
\end{pmatrix}.
$$
Let
$$
	\hat{\Psi}_k = QR\in \Re^{n\times l}
$$
be the \QR{} decomposition of $\hat{\Psi}_k$, where $Q\in \Re^{n \times n}$
has orthonormal columns and $R\in \Re^{n \times l}$
is upper triangular (see, e.g., \cite{GVL96}).

Then,
\begin{eqnarray*}
	B_{k+1} & = & B_0  + \Psi_k M_k \Psi_k^T  \\ &=&
	B_0 + \hat{\Psi}_kP^TM_kP \hat{\Psi}_k^T  \\
	&=& B_0 + Q R P^TM_kPR^TQ^T \\
\end{eqnarray*}

The matrix
$RP^TM_kPR^T$ is a real symmetric $n\times n$ matrix.
Since $\hat{\Psi}_k\in\Re^{n\times l}$, $R$ has at most rank $l$; moreover,
$R$ can be written in the form
$$R=\begin{pmatrix} R_1\\ 0 \end{pmatrix},$$
where $R_1\in\Re^{l\times l}$.
Then,
$$
RP^TM_kPR^T = \begin{pmatrix}R_1 \\ 0 \end{pmatrix} 
P^TM_kP
 \begin{pmatrix} R_1^T & 0 
\end{pmatrix}
= \begin{pmatrix}
R_1P^TM_kPR_1^T & 0 \\ 0 & 0
\end{pmatrix}.
$$
The eigenvalues of $RP^TM_kPR^T$ can be explicitly computed by forming the
spectral decomposition of $R_1P^TM_kPR_1^T\in\Re^{l\times l}$.  That is,
suppose $V_1D_1V_1$ is the spectral decomposition of $R_1P^TM_kPR_1^T$.
Then,
$$ 
RP^TM_kPR^T=
\begin{pmatrix}
R_1P^TM_kPR_1^T & 0 \\ 0 & 0
\end{pmatrix} 
= \begin{pmatrix}
V_1D_1V_1^T & 0 \\ 0 & 0
\end{pmatrix} = V D V^T
$$
where 
$$
	V \defined
	\begin{pmatrix}
		V_1 & 0 \\
		0 & I
	\end{pmatrix}\in\Re^{n\times n}
	\qquad \text{and} \qquad
	D \defined
	\begin{pmatrix}
		D_1 & 0 \\
		0 & 0
	\end{pmatrix}\in\Re^{n\times n}.
$$
This gives that
\begin{eqnarray}
	B_{k+1}
	&=& B_0 + Q V D V^T Q^T \nonumber \\
	&=& \gamma I + Q V D V^T Q^T \nonumber \\ 
	&=& Q V ( \gamma I+ D ) V^T Q^T, \label{eqn-eigs}
%	&=& B_0 + Q V D V^T Q^T \nonumber \\ 
%	&=& Q V ( B_0 + D ) V^T Q^T, \label{eqn-eigs}
\end{eqnarray}
yielding the spectral decomposition of $B_{k+1}$. 
% For $B_0=\gamma I$,
The matrix $B_{k+1}$ has an eigenvalue of $\gamma$ with multiplicity $n-l$ 
and $l$ eigenvalues given by
$\gamma + d_i$, where $1 \le i \le l$.
In practice, the matrices $Q$ and $V$ in (\ref{eqn-eigs}) are not stored.

In the case of the SR1{} updates, $\Psi_k=Y_k-B_0S_k$ and no permutation
matrix is required (i.e., $P=I$).

\medskip

Computing the eigenvalues of $B_{k+1}$ directly is an $O(n^3)$ process.  In
contrast, the above decomposition requires the \QR{} factorization of
$\Psi_k$ and the eigendecomposition of $R_1P^TM_kPR_1^T$, requiring
$(O(nl^2)$ flops) and $O(l^3)$ flops, respectively.  Since
$l\ll n$, the proposed method results in substantial computational savings.
%computing the eigenvalues of $B_{k+1}$ via
%QR factorization results in substantial computational savings.

%\todo{need: operation counts for everything--put in abstract and intro}\\
%\\

%\bigskip

%the thin qr + uniquess with r_ii>0 requires full rank

\medskip

The material presented above was first presented in~\cite{Burdakov13}
using the so-called ``thin QR'' factorization together with a Choleksy
factorization of an $m\times m$ symmetric matrix.
%Whereas, in the presentation above, the eigenvectors could be computed in
%theory by a matrix multiply.    
We chose to present the eigenvalue decomposition in terms of ``full QR''
since we are able update the \QR{} factorization as new 
quasi-Newton pairs are computed without having to
store $Q$ explicitly.  Here we describe in detail how the update to the
\QR{} factorization can be performed--enabling the eigenvalues of the
updated quasi-Newton matrix to be computed efficiently.

\subsection{Handling updates to $\hat{\Psi}$} \label{sec-updates}

In this section we detail handling updates to the \QR{} decomposition of
$\hat{\Psi}_k$ when additional limited-memory pairs are added to $S$ and
$Y$.  We consider two cases: Adding a limited-memory pair $(s_{k+1},y_{k+1})$ when
$k+1<m$ and when $k+1\ge m$, where $m$ is the maximum number of limited-memory
updates allowed to be stored.  The case $k+1<m$ requires adding a row and
column to the $R$ factor; whereas the case $k+1\ge m$ requires first deleting
a column (or two) of $\hat{\Psi}_k$ before adding the newest limited-memory
pair.  In both cases, the columns of $Q$ need not be formed nor stored.
However, when $\hat{\Psi}_k$ is not full rank, the \QR{} decomposition must
be computed from scratch.

We begin by discussing the process to compute $\hat{\Psi}_{k+1}$ from
$\hat{\Psi}_k$ when a new limited-memory pair is added to $S$ and $Y$.
The discussion considers the Broyden convex class of updates; however,
comments are included at the end of each subsection regarding the
\SR1{} case. 
%%%%%%%%%%%%%%%%%%%%%%%%%%%%%%%%%%%%%%%%%%%%%%%%%%%%%%%%%%%%%%%%%%%%
\subsubsection{Adding a column to $S$ and $Y$} \label{sec-add} Suppose
$\hat{\Psi}_k=QR\in\Re^{n\times l}$ is full rank and we have stored $k+1$
limited-memory Broyden convex class updates such that $k+1<m$, where $m$ is
the maximum number of limited-memory updates allowed to be stored by the
limited-memory quasi-Newton method.  Further, suppose we have computed a
$(k+2)$nd pair $(s_{k+1},y_{k+1})$.  To update the \QR{} decomposition, we
augment $\hat{\Psi}_k$ with the two columns $\begin{pmatrix} B_0 s_{k+1} &
  y_{k+1}\end{pmatrix}$.  This can be accomplished by using the procedure
proposed by Gill et al.~\cite{Gil74} for updating the \QR{} factorization
after a column is added.  This method relies upon $\hat{\Psi}_k$ having
full column rank.  For completeness, this procedure is presented below in
the context of adding two columns to $\hat{\Psi}_k$.  As in the previous
section, we assume that $B_k$ is updated using the Broyden convex set of
updates.

We begin by first adding the column $B_0s_{k+1}$ to $\hat{\Psi}_k$; the
same process may be followed to add the new last column, $y_{k+1}$.
Suppose 
\begin{equation} \label{eqn-psiP}
\hat{\Psi}_k = Q\begin{pmatrix} R_1 \\ 0 \end{pmatrix},
\end{equation}
where $R_1\in\Re^{l\times l}$.
Moreover, suppose we insert $B_0s_{k+1}$ into the final column of $\hat{\Psi}_k$
to obtain $\hat{\hat{\Psi}}_k$.
Setting $$\hat{\hat{\Psi}}_k = Q\begin{pmatrix} R_1 & u_1 \\ 0 & u_2\end{pmatrix}$$
yields that 
\begin{equation}\label{eqn-QR-add}
  B_0s_{k+1} = Qu \quad \text{with} \quad u=\begin{pmatrix} u_1 \\ u_2 \end{pmatrix},
\end{equation}
where $u_1\in\Re^{l}$ and $u_2\in\Re^{n-l}$.
  We now construct an orthogonal matrix $H_1$ such that 
\begin{equation}\label{eqn-QR-add2}
H_1
\begin{pmatrix} u_1 \\ u_2\end{pmatrix} =
\begin{pmatrix} u_1  \\ \eta \\ 0 \end{pmatrix},
\end{equation}
where $\eta=\pm \|u_2\|$, i.e., $$H_1 = \begin{pmatrix} I & 0 \\ 0 & \hat{H}_1
\end{pmatrix}$$ where $\hat{H}_1$ is a Householder matrix
such that $\hat{H}_1u_2 =\begin{pmatrix} \eta & 0 \end{pmatrix}^T$.
This choice of $H_1$ 
preserves the structure of $R_1$, i.e.,
$$
H_1\begin{pmatrix} R_1 & u_1 \\ 0 & u_2\end{pmatrix}
=
\begin{pmatrix} R_1 & u_1 \\ 0 & \eta \\ 0 & 0 \end{pmatrix}.
$$
Then, $\hat{\hat{\Psi}}_k=\hat{Q}\hat{R}$ is the \QR{} decomposition of $\hat{\hat{\Psi}}_k$,
where $$ \hat{R} = \begin{pmatrix} R_1 & u_1 \\ 0 & \eta \\ 0 & 0 \end{pmatrix},$$
and $\hat{Q}=QH_1^T$.

In this procedure, the matrices $Q$, $\hat{Q},$ and $H_1$ are not stored;
moreover, the unknowns $u_1$ and $\eta$ can be
computed without explicitly using these matrices.  Specifically, the relation in
(\ref{eqn-QR-add}) implies $\hat{\Psi}_k^TB_0s_{k+1} = \begin{pmatrix}R_1^T & 0 \end{pmatrix}Q^TQu$,
 i.e., 
\begin{equation}\label{eqn-u1}
\hat{\Psi}_k^TB_0s_{k+1} =R_1^Tu_1.
\end{equation} 
Equation (\ref{eqn-u1}) is a square $l\times l$ system that can be solved
for $u_1$ provided $\hat{\Psi}_k$ is full rank.  Finally, the scalar $\eta$
can be computed from the following relation obtained from combining
equations (\ref{eqn-QR-add}) and (\ref{eqn-QR-add2}): $\|B_0s_{k+1}\|^2
=\| \begin{pmatrix} u_1 & \eta\end{pmatrix}\|^2$. This yields that $\eta^2
= \|B_0s_{k+1}\|^2-\|u_1\|^2$.  This procedure can be repeated to add
$y_{k+1}$ to the new last column of $\hat{\hat{\Psi}}_k$, thereby updating
the \QR{} factorization of $\hat{\Psi}_k$ to $\hat{\Psi}_{k+1}$ with a new
pair of updates $(B_0s_{k+1}, y_{k+1})$.

The process of adding a new \SR1{} update to $\hat{\Psi}_k$ is simpler
since $\hat{\Psi}_k$ is augmented by only one column: $y_k-B_0 s_k$.

\subsubsection{The full-rank assumption}
The process described above requires $\hat{\Psi}_k$ to be full rank so that
there is a (unique) solution to (\ref{eqn-u1}).  When $\hat{\Psi}_k$ is not
full rank, the \QR{} decomposition must be computed from scratch.
Fortunately, there is an \emph{a priori} way to determine when there is no
unique solution: The matrix $\hat{\Psi}_k$ has full rank if and only if
$R_1$ in (\ref{eqn-u1}) is invertible; in particular, the diagonal of $R_1$
is nonzero.  When $R_1$ is singular, the process described above to update
the \QR{} decomposition for $\hat{\Psi}_{k+1}$ is skipped and the \QR{}
decomposition of $\hat{\Psi}_{k+1}$ should be computed from scratch at a
cost of  %about
$2l^2(n-l/3)$ flops. %$4kn(2k+1)$ flops \todo{check this!}.  
The process described
in Section~\ref{sec-add} can be reinstated to update the \QR{}
decomposition when the $R_1$ factor has nonzero diagonal entries, which may
occur again once the limited-memory updates exceed the maximum number
allowed, (i.e., $k\ge m$), and we are forced to delete the oldest pairs.

Similarly, when $\hat{\Psi}_k$ is ill-conditioned, $R_1$ will also be
ill-conditioned with at least one relatively small diagonal entry.  In this
case, (\ref{eqn-u1}) should not be solved; instead, the \QR{} factorization
should be computed from scratch.  As with the rank-deficient case, it is
possible to know this \emph{a priori}.

\subsubsection{Deleting and adding columns to $S$ and
  $Y$} \label{sec-fullrank} In this section, we detail the process to
update the \QR{} factorization in an efficient manner when $\hat{\Psi}_k$
is full rank and $k+1\ge m$.  As in the previous section, we assume we are
working with the Broyden convex class of updates.

Suppose $\hat{\Psi}_{k}=QR$ and we have stored the maximum number ($k+1$)
limited-memory pairs $\{(s_i,y_i)\}$, $i=0,\ldots, k$ allowed by the
limited-memory quasi-Newton method.  Further, suppose we have computed a
$(k+2)$nd pair $(s_{k+1}, y_{k+1})$. The process to obtain an updated \QR{}
factorization of $\hat{\Psi}_{k+1}$ from $\hat{\Psi}_k$ can be viewed as 
a two step process:
\begin{enumerate}
\item Delete a column of $S$ and $Y$.
\item Add a new column to $S$ and $Y$.
\end{enumerate}
%For each of these steps we track and update the \QR{} factorization.

\medskip

For the first step, we use ideas based on Daniel et al.~\cite{Dan76} and Gill et
al.~\cite{Gil74}.  Consider the Broyden class of updates. Suppose we
rewrite $\hat{\Psi}_{k}$ and $R$ as
\begin{equation}\label{eqn-Psi-delS}
\hat{\Psi}_k=
\begin{pmatrix} B_0 s_0 & y_0 & \tilde{\Psi}_k
\end{pmatrix}
\quad \text{and} \quad
R=\begin{pmatrix}
  r_1 & r_2 & \tilde{R}
\end{pmatrix},
\end{equation}
where $\tilde{\Psi}_k\in\Re^{n\times (l-2)}$ and
$\tilde{R} \in \Re^{n\times (l-2)}$.  This gives that 
$$
\hat{\Psi}_k =
\begin{pmatrix} B_0 s_0 & y_0 & \tilde{\Psi}_k
\end{pmatrix}
=Q\begin{pmatrix}
r_1 & r_2 & \tilde{R}
\end{pmatrix}.
$$
Deleting the first two columns of $\hat{\Psi}_k$ yields the matrix
$\bar{\Psi}_{k}=Q\tilde{R}$.
Notice that
$\tilde{R}$ has zeros beneath the second subdiagonal.
For clarity, we illustrate the nonzero entries of
$\tilde{R}$ for the case $n=8$ and $k=2$:
\begin{equation}\label{eqn-hatR}
\begin{pmatrix}
 * & * & * & * \\
 * & * & * & * \\
 * & * & * & * \\
 0 & * & * & * \\
 0 & 0 & * & * \\
 0 & 0 & 0 & * \\
 0 & 0 & 0 & 0 \\
 0 & 0 & 0 & 0 \\
\end{pmatrix},
\end{equation}
where * denotes possible nonzero entries.  
Givens rotations can be used to zero out the entries beneath the main
diagonal in $\tilde{R}$ at a cost of $24k^2-36k$.
(In the above example, eight entries
must be zeroed out to reduce (\ref{eqn-hatR}) to upper triangular form;
more generally, $4(k-1)$ entries must zeroed out to reduce $\tilde{R}$ to
upper triangular form.)
Let $G_{i,j}\in\Re^{n\times n}$ denote the Givens matrix that zeros out the
$(i,j)$th component of $\tilde{R}$, and suppose $\hat{G}$ is given by
$$\hat{G}\defined G_{2j-1,2j-2}G_{2j,2j-2}\cdots G_{2,1}G_{3,1}.$$
Then, $\hat{R}\defined\hat{G}\tilde{R}$ is an upper triangular matrix.
This yields the QR decomposition of the matrix
$\bar{\bar{\Psi}}_k$ defined as follows:
\begin{equation}\label{eqn-psihat-qr}
\bar{\bar{\Psi}}_{k}\defined\hat{Q}\hat{R},
\end{equation}
where $\hat{Q}=Q\hat{G}^T\in\Re^{n\times n}$ is orthogonal and
$\hat{R}\in\Re^{n\times (l-2)} $ is an upper triangular matrix.  With the
computation of $\bar{\bar{\Psi}}_k$ we have completed the first step.
Notice that neither $Q$ nor $\hat{Q}$ must be stored in order to
obtain $\hat{R}$.

\medskip

For the second step, the \QR{} factorization of
$\hat{\Psi}_{k+1}$ can be obtained from $\bar{\bar{\Psi}}_k$ using the
procedure outlined in Section~\ref{sec-add}.  \medskip

\medskip The process required for \SR1{} updates is simpler than for the
Broyden convex class of updates since it is only a rank-one update.  That is, 
only one column of $\hat{\Psi}_k$ must be deleted to remove the old
pair of updates and only one column must be added to incorporate the newest
pair of updates.

%%%%%%%%%%%%%%%%%%%%%%%%%%%%%%%%%%%%%%%%%%%%%%%%%%%%%%%%%%%%%%%%%%%%
\section{Numerical experiments}
In this section, we demonstrate the accuracy of the proposed method
implemented in \MATLAB{} to compute the eigenvalues of limited-memory
quasi-Newton matrices.  For the experiments, we considered limited-memory
\SR1{} matrices and three limited-memory members of the Broyden convex
class of updates; namely, limited-memory BFGS{} updates ($\phi=0$),
limited-memory \DFP{} updates ($\phi=1$), and limited-memory updates
obtained by selecting $\phi=0.5$ in (\ref{eqn-1param}).  The number of
limited-memory updates was set to 5 and $\gamma=3$.  The pairs $S$ and $Y$
were generated using random data.  In addition, in the case of the Broyden
convex class of updates, the updates were ensured to generate a positive
definite matrix by redefining $s_i$ as follows: $s_i=\text{sign}(s_i^Ty_i)s_i$
whenever $s_i^Ty_i<0$ for each $i\in\{1,\ldots,5\}$.  

We report the results of following three numerical experiments on each
quasi-Newton matrix:
\newline
\newline
\noindent
{\bf Experiment 1.} The eigenvalues of the quasi-Newton matrix were computed using the procedure outlined in Section 4.1.
\newline
{\bf Experiment 2.} A new quasi-Newton pair was generated, adding a column to both
$S$ and $Y$.  The procedure outlined in Section 4.2.1 was used
  to update $\Psi$, and then the eigenvalues were recomputed using the
  procedure outlined in Section 4.1.
\newline
{\bf Experiment 3.} The first columns of $S$ and $Y$ were deleted, simulating the
case when the oldest pair of updates is discarded.  Then, a new quasi-Newton
pair was generated, adding a new column to the end of both $S$ and $Y$. 
The procedure outlined in Section 4.2.3 was used to update $\Psi$, 
and then the eigenvalues were recomputed using the procedure
outlined in Section 4.1.

\medskip

To determine the accuracy of the proposed method, we explicitly formed each
quasi-Newton matrix and used the \MATLAB{} \texttt{eig} command to compute
its actual eigenvalues.  Due to memory
limitations in computing actual eigenvalues to test the proposed
method, we restricted the matrix sizes to $n\le 5000$.  
For each of the three
experiments, we report the size of the quasi-Newton matrix (``n'') and the
relative error in the computed eigenvalues measured by the infinity norm
(``RE Experiment 1'', ``RE Experiment 2'', and ``RE Experiment 3''); that
is, for each experiment the relative error was computed as
$$\text{RE} = \frac{\|(D+\gamma I) -
  \Lambda\|_{\infty}}{\|\Lambda \|_{\infty}},$$ where $(D+\gamma I)$ is as
in (\ref{eqn-eigs}) and $\Lambda$ is the matrix of eigenvalues obtained
using the \MATLAB{} \texttt{eig} command.

In Table 1, we report the results when $B$ was taken to be
randomly-generated limited-memory \SR1{} matrices of sizes $n=100, 500,
1000$ and $5000$.  For each matrix, the relative error in computing the
eigenvalues using the proposed method is very small (column 1).  As seen in
column 2, the relative error remained very small after a new column was
added to both $S$ and $Y$, $\Psi$ was updated using using Section 4.2.1,
and the eigenvalues were recomputed using the procedure outlined in Section
4.1.  Finally, column 3 shows that the relative error remained very small
after discarding the first stored quasi-Newton pair and adding a new column
to both $S$ and $Y$ using the procedures outlined in Section 4.2.3 and
Section 4.1.

\setlength\tabcolsep{1.5mm}
\begin{table}[!h]

  \caption{Summary of results when $B$ is a limited-memory SR1 matrix. } \begin{center}
  \begin{tabular}{|c|c|c|c|}
    \hline $n$ & RE Experiment 1 &
    RE Experiment 2  & RE Experiment 3\\
\hline \texttt{100} & \texttt{1.92439e-15} & \texttt{2.07242e-15} & \texttt{2.81256e-15} \\
\hline \texttt{500} & \texttt{4.88498e-15} & \texttt{4.44089e-15} & \texttt{6.21725e-15}\\
\hline \texttt{1000} & \texttt{8.14164e-15} & \texttt{7.99361e-15} & \texttt{7.84558e-15}\\
\hline \texttt{5000} & \texttt{1.71714e-14} & \texttt{1.98360e-14} & \texttt{1.68754e-14}
\\
\hline
\end{tabular} 
\end{center}
\end{table}

Table 2 reports the results when $B$ was a randomly-generated limited-memory
\BFGS{} matrix of sizes $n=100, 500, 1000,$ and $5000$.  In all cases,
the proposed method computed the eigenvalues to high accuracy.
Table 3 reports the results when $B$ was a randomly-generated
limited-memory \DFP{} matrix of various sizes.  As in Tables 1 and 2,
the proposed method computed the eigenvalues of these matrices
to high accuracy in each experiment.

\setlength\tabcolsep{1.5mm}
\begin{table}[!h]
\caption{Summary of results when $B$ is a limited-memory BFGS matrix.}
\begin{center}
			\begin{tabular}{|c|c|c|c|}
    \hline $n$ & RE Experiment 1 &
    RE Experiment 2  & RE Experiment 3\\
\hline \texttt{100} & \texttt{5.53332e-16} & \texttt{1.21039e-16} & \texttt{7.86896e-16} \\
\hline \texttt{500} & \texttt{6.35220e-16} & \texttt{4.28038e-16} & \texttt{5.86555e-16}\\
\hline \texttt{1000} & \texttt{1.13708e-15} & \texttt{2.39590e-15} & \texttt{1.62325e-15}\\
\hline \texttt{5000} & \texttt{1.14773e-15} & \texttt{3.39882e-15} & \texttt{1.30101e-15}\\
\hline
                      \end{tabular} 
\end{center}
\end{table}

\setlength\tabcolsep{1.5mm}
\begin{table}[!h]
  
  \caption{Summary of results when $B$ is a limited-memory DFP matrix.}
  \begin{center}
			\begin{tabular}{|c|c|c|c|}
    \hline $n$ & RE Experiment 1 &
    RE Experiment 2  & RE Experiment 3\\
\hline \texttt{100} & \texttt{1.69275e-15} & \texttt{2.05758e-16} & \texttt{3.65114e-16} \\
\hline \texttt{500} & \texttt{9.58309e-16} & \texttt{6.19241e-16} & \texttt{2.10460e-15}\\
\hline \texttt{1000} & \texttt{4.15522e-15} & \texttt{1.30844e-14} & \texttt{1.72417e-14}\\
\hline \texttt{5000} & \texttt{2.27937e-15} & \texttt{1.20206e-14} & \texttt{2.97026e-15}
\\
\hline
                      \end{tabular} 
\end{center}
\end{table}

Finally, Table 4 reports the results when $B$ obtained using the Broyden
convex class of updates with $\phi=0.5$.  In all experiments with this type
of update, the proposed method was able to compute all the eigenvalues to
high accuracy.

\setlength\tabcolsep{1.5mm}
\begin{table}[!h]

  \caption{Summary of results when $B$ is a limited-memory member of the Broyden class of convex updates with $\phi=0.5$.}
\begin{center}
			\begin{tabular}{|c|c|c|c|}
    \hline $n$ & RE Experiment 1 &
    RE Experiment 2  & RE Experiment 3\\
\hline \texttt{100} & \texttt{5.11757e-15} & \texttt{9.05737e-15} & \texttt{6.02940e-16}\\
\hline \texttt{500} & \texttt{1.11222e-15} & \texttt{4.90513e-15} & \texttt{1.60814e-15}\\
\hline \texttt{1000} & \texttt{1.76830e-15} & \texttt{2.83112e-15} & \texttt{2.18559e-15}\\
\hline \texttt{5000} & \texttt{9.86622e-15} & \texttt{2.95003e-15} & \texttt{5.88569e-15} 
\\
\hline
                      \end{tabular} 
\end{center}
\end{table}

%%%%%%%%%%%%%%%%%%%%%%%%%%%%%%%%%%%%%%%%%%%%%%%%%%%%%%%%%%%%%%%%%%%%
\section{Concluding remarks}

In this paper we produced the compact formulation of quasi-Newton matrices
generated by the Broyden convex class of updates. Together with the \QR{}
factorization, this compact representation was used to compute the
eigenvalues of any member of this class of updates.  In addition, we
presented an efficient procedure to update the QR factorization when a new
pair of updates for the quasi-Newton matrix is computed.  With this
approach we are able to substantially reduce the computational costs of
computing the eigenvalues of quasi-Newton matrices.  Applications of this
work are the subject of current research.  Code and drivers used
for this paper can be found on the following website:
$$\texttt{http://users.wfu.edu/erwayjb/software.html}.$$

%%%%%%%%%%%%%%%%%%%%%%%%%%%%%%%%%%%%%%%%%%%%%%%%%%%%%%%%%%%%%%%%%%%%

\renewcommand{\appendix}{%
  \par
  \setcounter{section}{0}%
  \renewcommand{\thesection}{A\Alph{section}}%
}

\appendix
\include{Appendix  \textbf{Appendix A.}}
\numberwithin{equation}{section}
In~\cite{ApoSP11, ApoSP08}, Apostolopoulou et al. find explicit formulas
for computing the eigenvalues of a \BFGS{} matrix when {\bf at most} two
limited-memory quasi-Newton pairs are used to update an initial $B_0$.
While the methods in this paper are not limited to two updates and can be
applied to \SR1{} matrices and any matrix generated using the Broyden
convex class of updates, we show that the results found
in~\cite{ApoSP11, ApoSP08} for the case of one update can be derived using
the technique proposed in this paper.

Without loss of generality, Apostolopoulou et al. derive a formula
for computing the eigenvalues of the following matrix obtained after applying one update:
\begin{equation}\label{eqn-A1}
	B_1 = \frac{1}{\theta_{0}} I - \frac{1}{\theta_{0}} \frac{s_0s_0^T}{s_0^Ts_0} +
	\frac{y_0y_0^T}{s_0^Ty_0}, \quad \text{where \ }
	\theta_{0} = \frac{s_0^Ts_0}{s_0^Ty_0}.
\end{equation} 

The compact formulation of $B_1$ is given by $B_1=B_0+\Psi_0M_0\Psi_0^T$
where 
$$\Psi_1 = 
	\begin{bmatrix}
		B_0 s_0 & y_0
	\end{bmatrix} \quad \text{and} \quad
M_0=	\begin{bmatrix}
		-s_0^TB_0 s_0 & 0 \\
		0 & s_0^Ty_0
	\end{bmatrix}^{-1}.
$$
For notational simplicity, we drop the subscript $k=0$ for the duration of
Appendix A.  We now compute the eigenvalues of $B_1 = (1/\theta)I + \Psi M
\Psi^T$ using the \QR{} factorization of $\Psi$ and show that these
eigenvalues are the same as those obtained in~\cite{ApoSP11, ApoSP08}
by comparing the characteristic polynomials.

The QR factorization can be computed using Householder transformations.
The first Householder transformation zeros out all the elements in the
first column of $\Psi$
below the first entry.  
Let $v_1 \defined  \Psi e_1 -\|\Psi e_1\|e_1  $, where $e_1$ denotes the
first canonical basis vector.  Since $\Psi e_1 = (1/\theta)s$, 
then 
$$Q_1 = I - \frac{2}{v_1^Tv_1} v_1v_1^T$$ is such that
$$ Q_1 \left(\Psi e_1\right) = Q_1\left(\frac{1}{\theta}\right)s = \frac{\|s\|}{\theta} e_1.$$
(For more details on constructing Householder matrices, see e.g.,~\cite{GVL96}.)

Using the definition of $Q_1$ and $v_1$ gives that
\begin{equation}\label{eqn-hh1}
Q_1 \Psi = Q_1 [ (1/\theta)s \ \ \ y] = \begin{bmatrix} \frac{\|s\|}{\theta}e_1 \ \ \ y-\frac{2v_1^Ty}{v_1^Tv_1} v_1\end{bmatrix} = 
\left[\begin{array}{c|c}
\frac{\|s\|}{\theta} &  \\
 0                   &  y-\frac{2v_1^Ty}{v_1^Tv_1} v_1 \\ 
\vdots               &  \\
0                    &  
\end{array} \right]. \end{equation}
%		\ \frac{\|s\|}{\theta}  \ & y^{(1)}-  \left ( \frac{y^Ts - \| s \|y^Te_1 }{s^Ts - \| s \| s ^Te_1}  \right ) ( s^{(1)} - \| s\|)\\
%		0 & y^{(2)} - \left ( \frac{y^Ts - \| s \|y^Te_1 }{s^Ts - \| s \| s ^Te_1}  \right )  s^{(2)}\\
%		\vdots & \vdots \\
%		0 & y^{(n)} -   \left ( \frac{y^Ts - \| s \|y^Te_1 }{s^Ts - \| s \| s ^Te_1}  \right )   s^{(n)}\\
%	\end{bmatrix}.

The second Householder transformation zeros out all entries in the second
column below the second row. Let
$$u_2 = 
\begin{bmatrix}
(y-\frac{2v_1^Ty}{v_1^Tv_1} v_1)^Te_2\\
\vdots\\
(y-\frac{2v_1^Ty}{v_1^Tv_1} v_1)^Te_n
\end{bmatrix}.
$$
Defining $v_2 \defined u_2 - \|u_2\|e_1$, 
then $$\hat{Q}_2 = I-\frac{2}{v_2^Tv_2}v_2v_2^T$$
is such that $\hat{Q}_2 u_2 = \|u_2\|e_1$, and
$$Q_2 \defined \begin{bmatrix} 1 & 0 \\ 0 & \hat{Q}_2\end{bmatrix}$$
is such that
\begin{equation}\label{eqn-hh2}
R\defined	Q_2Q_1 \Psi = 
	\begin{bmatrix}
\frac{\|s\|}{\theta} &  \left(y-\frac{2v_1^Ty}{v_1^Tv_1} v_1\right)^Te_1\\
 0                   &  \|u_2\| \\ 
 0                   &    0 \\
\vdots               &   \vdots        \\
0                    &  0 
\end{bmatrix}
\end{equation}
is an upper triangular matrix.
If $Q\defined Q_1^TQ_2^T$ then by (\ref{eqn-hh2})
\begin{equation}\label{eqn-BcompactA}
	B_1  = Q \left ( \frac{1}{\theta} I + RMR^T \right ) Q^T.
\end{equation}
Let $R_1 \in \Re^{2 \times 2}$ be the first $2 \times 2$ block of $R$, i.e.,
$$
	R_1 = 
	\begin{bmatrix}
\frac{\|s\|}{\theta} &  \left(y-\frac{2v_1^Ty}{v_1^Tv_1} v_1\right)^Te_1\\
 0                   &  \|u_2\|
	\end{bmatrix}. 
$$
The matrix $R_1$ can be further simplified by noting that
\begin{eqnarray*}
\left(y-\frac{2v_1^Ty}{v_1^Tv_1} v_1\right)^Te_1
& = & 
\frac{y^Te_1v_1^Tv_1}{v_1^Tv_1}-\frac{2v_1^Ty}{v_1^Tv_1} v_1^Te_1 \\
& = & 
\frac{y^Te_1(s-\|s\|e_1)^T(s-\|s\|e_1)  -2(s-\|s\|e_1)^Ty(s-\|s\|e_1)^Te_1}{(s-\|s\|e_1)^T(s-\|s\|e_1)}  \\
& = & 
\frac{2y^Te_1\|s\|(\|s\|-s^Te_1)+2(s-\|s\|e_1)^Ty(\|s\|-s^Te_1)}{2\|s\|(\|s\|-s^Te_1)} \\
& = & 
\frac{y^Te_1\|s\|+(s-\|s\|e_1)^Ty}{\|s\|} \\
& = & 
\frac{s^Ty}{\|s\|},
\end{eqnarray*}
and thus,
$$
	R_1 = 
	\begin{bmatrix}
\frac{\|s\|}{\theta} &   \frac{s^Ty}{\|s\|}.\\ 
 0                   &  \|u_2\|
	\end{bmatrix}. 
$$

Rewriting (\ref{eqn-BcompactA}) using $R_1$ yields
$$
	B_1  = Q
	\begin{bmatrix}
		\frac{1}{\theta} I_2 + R_1MR_1^T & 0 \\
		0 & \frac{1}{\theta}I_{n-2}
	\end{bmatrix}
	Q^T,
$$
implying that the eigenvalues of $B_1$ are the union of eigenvalues of the
two diagonal blocks in $B_1$.  Note that the leading $2\times 2$ block
can be simplified as follows:
\begin{eqnarray*}
 \frac{1}{\theta} I_2 + R_1MR_1^T  & = & 
	\frac{1}{\theta} I -
	\begin{bmatrix}
		 \frac{\|s\|}{\theta}   & \frac{s^Ty}{\|s\|}\\
		0 &  \|u_2\|
	\end{bmatrix}
	\begin{bmatrix}
		\frac{\theta}{s^Ts} & 0 \\
		0 & -\frac{1}{s^Ty}
	\end{bmatrix}
	\begin{bmatrix}
		\ \frac{\|s\|}{\theta}  \ & 0 \\
\frac{s^Ty}{\|s\|}  &  \|u_2\|
	\end{bmatrix}\\ 
	&=& \frac{1}{\theta}I - 
	\begin{bmatrix}
		\frac{1}{\theta}  -  \frac{(s^Ty)}{\|s\|^2} &  -   \frac{\|u_2\|}{\|s\|} \\
-  \frac{\|u_2\|}{\|s\|} 
& -\frac{\|u_2\|^2}{s^Ty}
	\end{bmatrix}\\
	&=& 
	\begin{bmatrix}
	 \frac{(s^Ty)}{\|s\|^2} &  \frac{\|u_2\|}{\|s\|} \\ 
 \frac{\|u_2\|}{\|s\|}
& \frac{1}{\theta} + \frac{\|u_2\|^2}{s^Ty}
	\end{bmatrix}. \label{eqn-2x2}
\end{eqnarray*}

The characteristic polynomial of leading $2\times 2$ block of
$B_1$ is given by
\begin{equation}\label{eqn-QQ}
	\text{det} \left ( \frac{1}{\theta} I + R_1 M R_1^T - \lambda I \right )
= \lambda^2 - \lambda\left(\frac{1}{\theta} + \frac{\|u_2\|^2}{s^Ty} +
\frac{s^Ty}{\|s\|^2}
\right) %end of lambda
+ \frac{1}{\theta}\frac{s^Ty}{\|s\|^2}.
\end{equation}
Finally, since 
\begin{equation*}
	\left[\left(y-\frac{2v_1^Ty}{v_1^Tv_1} v_1\right)^Te_1\right]^2 + \|u_2\|^2	 =\| Q_2Q_1 y \|^2  =		 y^Ty \quad \text{and}
\quad \left[\left(y-\frac{2v_1^Ty}{v_1^Tv_1} v_1\right)^Te_1\right]^2 = \frac{(s^Ty)^2}{\|s\|^2},
\end{equation*}
then (\ref{eqn-QQ}) simplifies to
\begin{equation}\label{eqn-QQ2}
\lambda^2-\frac{\lambda}{\theta}\left(1 + \theta\frac{y^Ty}{s^Ty} %end of lambda
+\frac{1}{\theta^2}
\right).
\end{equation}
Thus, the characteristic polynomial of $B_{1}$ is given by
$$
	p(\lambda) = \left(\lambda^2-\frac{\lambda}{\theta}\left(1 + \theta\frac{y^Ty}{s^Ty}\right)
	+ \frac{1}{\theta^2} \right)\left( \lambda - \frac{1}{\theta} \right )^{n-2},
$$
which is the same as the characteristic polynomial derived in~\cite[Equation 4]{ApoSP11} and 
~\cite[Equation 9]{ApoSP08}.

%%%%%%%%%%%%%%%%%%%%%%%%%%%%%%%%%%%%%%%%%%%%%%%%%%%%%%%%%%%%%%%%%%%%
%\section*{References}

\bibliographystyle{abbrv}
\bibliography{ErwayMarcia}

\end{document}